\documentclass[10pt,twoside]{article}
\usepackage{graphicx}
\usepackage{amsfonts}
\usepackage{amsmath}
\usepackage{Latex-document}

\title{\bf  List Colouring of Graphs with at Most\vskip -1mm {\boldmath
$\big(2-o(1)\big)\chi$} Vertices \vskip 6mm}
\author{{\bf Bruce Reed}\thanks{CNRS, Paris, France and School
of Computer Science, McGill University, Montreal, Canada.E-mail: breed@jeff.cs.mcgill.ca. This research was
partially supported by DIMACS and by a CNRS/NSF collaboration grant.} \quad Benny Sudakov
\thanks{Department of Mathematics,
 Princeton University, Princeton, NJ 08540, USA
and Institute for Advanced Study, Princeton, NJ 08540, USA. E-mail: bsudakov@math.princeton.edu. Research
supported in part by NSF grants DMS-0106589, CCR-9987845 and by the State of New Jersey.}\vspace*{-0.5cm}}
\date{\vspace{-8mm}}
\newtheorem{theo}{Theorem}

\newtheorem{observation}[theo]{Observation}
\newtheorem{lemma}[theo]{Lemma}

\newcommand{\proofend}{\hspace*{\fill}\mbox{$\Box$}\vspace{2ex}}

\begin{document}
\maketitle

\thispagestyle{first} \setcounter{page}{587}

\begin{abstract}\vskip 3mm
Ohba has conjectured \cite{ohb} that if the graph $G$ has
$2\chi(G)+1$ or fewer vertices then the list chromatic number and
chromatic number of $G$ are equal. In this paper we prove that
this conjecture is asymptotically correct. More precisely we
obtain that for any $0<\epsilon<1$, there exist an
$n_0=n_0(\epsilon)$ such that the list chromatic number of $G$
equals its chromatic number,  provided
$$n_0\,\leq\,|V(G) |\, \le\, (2-\epsilon)\chi(G).$$
\vskip 4.5mm


\end{abstract}

\vskip 12mm

\section{Introduction} \label{section 1}\setzero
\vskip-5mm \hspace{5mm}

Recently, a host of important results on graph colouring have been
obtained via the probabilistic method. The first author presented
an invited lecture at the 2002 International Congress of
Mathematicians surveying a number of these results. The recent
monograph \cite{molreed}  provides a more in depth survey of the
topic. This paper presents one example of a result proven using
the method.

An instance of List Colouring  consists of a graph $G$ and a list
$L(v)$ of colours for each vertex $v$ of $G$. We are asked to
determine if there is an {\it acceptable}  colouring of $G$, that
is a colouring in which each vertex receives a colour from its
list, and no edge has both its endpoints coloured with the same
colour. The {\em list-chromatic number} of $G$, denoted
$\chi^l(G)$  is the minimum integer $k$ such that for every
assignment of a list $L(v)$ of size at least $k$ to every vertex
$v$ of $G$, there exist an acceptable colouring of $G$. The
list-chromatic number was introduced by Vizing \cite{viz}, and
independently Erd\"{o}s et al. \cite{erd}. This parameter has
received a considerable amount of attention in recent years (see,
e.g. \cite{jen}, \cite{alo}).

\begin{figure}
\begin{center}
\includegraphics
[width=.5\textwidth,height=.22\textheight] {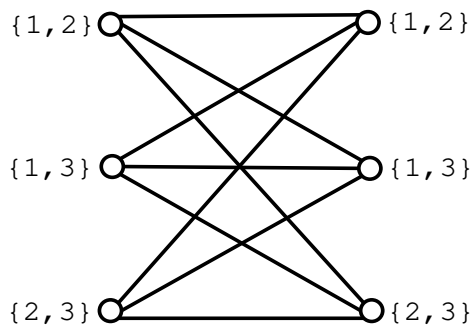} \caption[]{A bipartite graph with list chromatic number
three} \label{fst}
\end{center}
\end{figure}

Clearly, by definition, $\chi^l (G) \ge \chi(G)$ because
$\chi(G)=k$ precisely if an acceptable colouring exists when each
$L_v$ is $\{1,...,k\}$. However, the converse inequality is not
true, e.g.  $\chi^l (K_{3,3})=3$ as can be easily verified by
considering Figure \ref{fst}. In fact, there are bipartite graphs
with arbitrarily high chromatic number (indeed even for bipartite
$G$, $\chi^l(G)$ is bounded from below by a function of the
minimum degree which goes to infinity, see \cite{alo}). This shows
that the gap between $ \chi(G)$ and $\chi^l (G)$ can be
arbitrarily large. Moreover it shows that $\chi^l (G)$ can not be
bound by any function of the chromatic number of $G$. This gives
rise to the following intriguing question in the theory of graph
colourings: Find conditions which guarantee the equality of the
chromatic and list-chromatic numbers.

There are many conjectures hypothesizing conditions on $G$ which
imply that $\chi(G)=\chi^l(G)$. Probably, the most famous of these
is the List Colouring Conjecture (see \cite{jen}) which states
that this is true if $G$ is a line graph. One interesting example
of a graph with $\chi=\chi^l$ was obtained in the original paper
of Erd\"{o}s et al. \cite{erd}. They proved that if $G$ is
complete $k$-partite graph with each part of size two then
$\chi(G)=\chi^l(G)=k$. It took nearly twenty years until Ohba
\cite{ohb} noticed that this example is actually part of much
larger phenomenon. He conjectured (cf. \cite{ohb}) that $\chi(G) =
\chi^l (G)$ provided $|V(G) | \le 2\chi(G)+1$. This conjecture if
it is correct is best possible. Indeed, let $G$ be a complete
$k$-partite graph with $k-1$ parts of size $2$ and one part of
size $4$. Then the number of vertices of $G$ is $2k+2$, the
chromatic number is $k$ and it was proved in \cite{ohb} that
list-chromatic number of $G$ is at least $k+1>k$.

In his original paper Ohba obtained that $\chi^l(G) = \chi(G)$
for all graphs $G$ with $|V(G) | \le \chi(G)+\sqrt{2\chi(G)}$. His
conjecture was settled for some other special cases in \cite{eno}.
Recently the result of Ohba was substantially improved by the
authors of this paper. In \cite{RS} they proved that Ohba's
conjecture is true for all graphs $G$ with at most ${5 \over
3}\chi(G)-{4 \over 3}$ vertices. In this paper we want to improve
this result for large graphs and prove that the  conjecture is
asymptotically correct. More precisely we obtain the following
theorem.

\begin{theo}
\label{maint} For any $0<\epsilon<1$, there exist an
$n_0=n_0(\epsilon)$ such that $\chi^l(G) = \chi(G)$ provided $n_0
\leq |V(G) | \le (2-\epsilon)\chi(G)$.
\end{theo}

The rest of this paper is organized as follows. In the next
section we describe the main steps in the proof of Theorem
\ref{maint}. More precisely, we present our key lemma and show how
to deduce from it the assertion of the theorem. We will prove this
lemma using probabilistic arguments. In Sections 3 and 4 we
discuss the main ideas we are going to use in the proof. We
present the details of the proof in Section 5. Finally, the last
section of the paper contains some concluding remarks.

\section{The key lemma}
\label{section 2} \setzero\vskip-5mm \hspace{5mm}

In this section we present the main steps in the proof of Theorem
\ref{maint}. First we need the following lemma from \cite{RS}
whose short proof we include here for the sake of completeness.

\begin{lemma}
\label{lessn} For any integer $t$, if $\chi^l(G)>t$ then there
exist a set of lists $L(v), v \in V(G)$ for which there is no
acceptable colouring such that each list has at least $t$ elements
and the set ${\cal A} = \cup_{v \in V(G)} L(v) $ has size less
than $|V(G)|$.
\end{lemma}

\noindent {\bf Proof.}\, Assume $\chi^l(G) >t$ and choose a set of
lists $L(v), v \in V(G)$ for which there is no acceptable
colouring, in which each list has size at least $t$ and which
minimizes $|{\cal A}|$.

Now, if $|{\cal A}| < |V(G)|$ then we are done. So, we can assume
the contrary. We consider the bipartite graph $H$ with bipartition
$({\cal A}, V(G))$ and an edge between $c$ and $v$ precisely if $c
\in L(v)$. If there is a matching of size $|V|$ in $H$ then this
matching saturates $V$ and points out an acceptable colouring for
the List Colouring instance in which no colour is used more than
once. Since, there is no such acceptable colouring, no such
matching exists. Thus there must be a smallest subset $B$ of
${\cal A}$ which is not the set of endpoints of a matching in this
graph and this set must have at most $|V|$ elements. Clearly, $B$
contains at least two vertices. Now, by the minimality of $B$
there is a matching $M$ in $H$  of size $|B|-1$ whose endpoints in
${\cal A}$ are in $B$. Further, classical results in matching
theory (see e.g. Theorem 1.1.3 of \cite{lp}) tell us that if $W$
is the set of endpoints of $M$ in $V$ then for $v \not \in W$, we
have $L(v) \cap B = \emptyset$.

Let $x$ be any vertex in $G-W$ and replace $L(v)$ by $L(x)$ for
every vertex $v\in W$. This yields a new List Colouring Problem in
which the total number of colours in  all lists is smaller than
$|{\cal A}|$ (since all the new list are disjoint from B).
Therefore by the minimality of our original choice, there exist an
acceptable colouring of $G$ for this new Lists Colouring instance.
In particular this implies that we can obtain an acceptable
colouring of $G-W$ for the original lists $L(v)$. Since no colour
in $B$ is used in this colouring, using the colouring of $W$
pointed out by $M$ yields an extension of this colouring to a
colouring of $G$ in which no colour of $B$ appears more than once.
This contradicts our assumption that there is no acceptable
colouring for this instance and proves the lemma. \proofend

\noindent {\bf Proof of Theorem \ref{maint}.}\, Let $0<\epsilon<1$
be a fixed constant and let $G$ be a graph satisfying
$|V(G)|<(2-\epsilon)\chi(G)$. We assume that $\chi^l(G)>\chi(G)$
and obtain a contradiction. Since adding an edge between vertices
in different colour classes in an optimal colouring of $G$ does
not change $\chi(G)$ and can only increase $\chi^l(G)$, we will
assume that $G$ is complete $\chi(G)$-partite graph. Thus $G$ has
a unique partition into $\chi(G)$ stable sets. We refer to these
stable sets as parts rather than colour classes so as to avoid
confusion with the colours used in our acceptable colouring of
$G$.

Now by Lemma \ref{lessn}, if $\chi^l(G)> \chi(G)$ then there is an
instance of List Colouring  on $G$ for which no acceptable
colouring exists, in which each list has length at least $\chi(G)$
and  such that the size of the union of all lists $L(v)$ is less
than $|V(G)|$. This means that in an acceptable colouring at least
one colour must be used  on more than one vertex. Fortunately, it
also implies that for every non-singleton part $U$ there is at
least one colour which appears on $L(v)$ for more than half the
vertices of $U$ (since each $L(v)$ contains more than half the
colours).

Our proof approach is simple. For each non-singleton part $U$, we
choose some colour $c_U$ and colour with $c_U$ all the vertices of
$U$ whose list contains $c_U$ (thus we must insist that all the
$c_U$ are distinct). We  complete the colouring by finding a
bijection between the vertices not yet coloured and the colours
not yet used so that each such colour is in the list of the vertex
with which it is matched. This yields an acceptable colouring in
which for each part $U$ there is at most one colour $c_U$ used on
more than one vertex of $U$.

To begin, we consider the case when there is some part $U$ such
that some colour appears on all the vertices of $U$. We show that
we can reduce to a smaller problem by using any such colour for
$c_U$. Iteratively repeating this process yields a graph where no
such $U$ exists and hence, in particular, there are no parts of
size two.

Our choices for the remaining $c_U$ are discussed in the proof of
the key Lemma \ref{bruce} which consists of the analysis of a
probabilistic procedure for choosing the remaining $c_U$.
Unfortunately, before discussing this procedure we need to deal
with some technical details.

So, to begin we show that we can assume that
 $\cap_{v\in U}L(v)$ is empty for all
parts $U$ of size bigger than $1$ in the partition of $G$. To see
this, let $U$ be a part of size at least $2$ such that $\cap_{v\in
U}L(v)\not = \emptyset$. Then  the graph $G-U$ has chromatic
number $\chi(G)-1$ and  at most $|V(G)|-2$ vertices  and therefore
also satisfies
$$|V(G-U)| \leq |V(G)|-2 \leq (2-\epsilon)\chi(G)-2=(2-\epsilon)(\chi(G)-1)-\epsilon <(2-\epsilon)\chi(G-U).$$
Note that it also satisfies $\chi^l(G-U)> \chi(G-U)=\chi(G)-1$
since otherwise we can obtain an acceptable colouring of $G$ from
the lists $L(v)$. Indeed, let $c$ be a colour in $\cap_{v\in
U}L(v)$. Since $\chi^l(G-U) =\chi(G)-1$, we know there is an
acceptable colouring of $G-U$ from the lists $L(v)-c$. Colouring
all vertices in $U$ with $c$ we obtain an extension of this
colouring to an acceptable colouring of $G$  from the original
lists,  a contradiction. Therefore we will consider the graph
$G-U$ instead of $G$ and continue this process until we obtain a
graph $G'$ and an instance of List Colouring  on $G'$ with the
following properties.
\begin{itemize}
\item
$G'$ is $\chi(G')$-partite graph which satisfies
$|V(G')|<(2-\epsilon)\chi(G')$.
\item
Each list $L'(v)$ has length at least $\chi(G')$ and there is no
acceptable colouring of $G'$ from $L'(v)$.
\item
The size of the union of all lists is less than $|V(G')|$.
\item
$\cap_{v\in U}L'(v)$ is empty for all parts of size bigger than
$1$ in the partition of $G'$.
\end{itemize}

Since the size of the lists is $\chi(G')>|V(G')|/2$ we obtain that
$L'(x) \cap L'(y)$ is non empty for any two vertices $\{x, y\}$ in
$G'$. In particular this implies that in the partition of $G'$
there are no parts of size two. Note that the original graph $G$
has at most $|V(G)|/2$ parts of size $\geq 2$ and each time we
removed such  a part the chromatic number of the remaining graph
decreased by one. Therefore we decrease chromatic number of $G$ by
at most $|V(G)|/2$ and hence the remaining graph $G'$ should have
at least
$$|V(G')| \geq \chi(G') \geq \chi(G)-\frac{|V(G)|}{2} \geq \frac{|V(G)|}{2-\epsilon}-\frac{|V(G)|}{2}\geq\frac{\epsilon|V(G)|}{4}$$
vertices. So by choosing an appropriate bound on the size of
$|V(G)|$ we can make $|V(G')|$ arbitrarily large. This completes
our discussion of parts $U$ for which some colour is in $L(v)$ for
all vertices $v$ of $U$. We turn now to the technical details
necessary before we present the rest of the ideas needed in the
proof.

Let $X$ be the  set of all the vertices in the singleton classes
in the partition of $G'$. Pick $m$ to be a sufficiently large
integer constant $m=m(\epsilon)$ and let $t$ be an integer which
satisfies

\begin{eqnarray}
\label{A} \frac{t+1}{m} \leq \frac{|X|}{\chi(G')} \leq
\frac{t+2}{m}.
\end{eqnarray}
Since in the partition of $G'$ there are no parts of size two, we
obtain that $|X|+3(\chi(G')-|X|) \leq
|V(G')|<(2-\epsilon)\chi(G')$. This implies that $|X|
\geq(1+\epsilon)\chi(G')/2$ and that $m/2<t\leq m-2$.

Set ${\cal A}=\cup_{v\in G'}L'(v)$. Let $H$ be  a bipartite graph
with bipartition $(X,{\cal A})$ and an edge between $c$ and $v$
precisely when $c \in L'(v)$. Note that the degree of every vertex
from $X$ in $H$ is at least $\chi(G')\geq |V(G')|/2>|{\cal A}|/2$.
Therefore by well known results on Zarankiewicz's problem (see,
e.g., \cite{Lo}, Problem 10.37), $H$ contains a complete bipartite
graph with $t$ vertices in $X$ and $m$ vertices in $\cal A$.
Denote the set of vertices from $X$ and $\cal A$ by ${\cal S}_1$
and ${\cal C}_1$ respectively and remove them from $H$. Note that
the bound on Zarankiewicz's problem guarantees that we will
continue to find a copy of the complete bipartite graph $K_{t,m}$
in $H$ until the minimal degree of a vertex in $X$ is
$o(\chi(G'))=o(|{\cal A}|).$ Thus in the end we obtain at least
$k=(1-o(1))\chi(G')/m $ disjoint sets of colours ${\cal C}_1,
\ldots, {\cal C}_k$ and also $k$ disjoint sets of singleton
partition classes ${\cal S}_1, \ldots, {\cal S}_k$, such that
${\cal C}_i \subset L'(s)$ for every vertex $s \in {\cal S}_i$.
Denote by ${\cal C}=\cup_i {\cal C}_i$, by  ${\cal S}=\cup_i {\cal
S}_i$ and let $C$ and $S$ be the sizes of ${\cal C}$ and ${\cal
S}$ respectively. Now using  (\ref{A}) we can obtain the following
inequalities
\begin{eqnarray*}
|X|- S&=&|X|-kt \leq
\frac{t+2}{m}\chi(G')-(1+o(1))\frac{t}{m}\chi(G')\\
&=& (1+o(1))\frac{2}{m}\chi(G')=(2+o(1))k<3k
\end{eqnarray*}
and
$$ |X|-S=|X|-kt \geq \frac{t+1}{m}\chi(G')-(1+o(1))\frac{t}{m}\chi(G')=(1+o(1))\frac{\chi(G')}{m}.$$
In the above discussion and in particular in the last two
inequalities we used that $m$ and $t$ are constants but $|V(G')|$
(and thus also $\chi(G')$) tends to infinity.

Let $W$ be the union of some set of $r=\chi(G')-|{\cal C}|$
singleton partition classes which do not belong to $\cal S$.
Such  a set  $W$ exists, since  the number of singleton partition
classes outside $\cal S$ is at least $(1+o(1))\chi(G')/m \gg r
=\chi(G')-km=o(\chi(G'))$. Note that we can obtain an acceptable
colouring of $W$ with the lists $L'(v)-{\cal C}$ greedily, since
the size of $L'(v)-{\cal C}$ is equal to $r$.  Let  $T$ be the set
of $r$ colours used to colour $W$ in one such acceptable
colouring. Denote by $G''=G'-W$ and let $L''(v)=L'(v)-T$ for every
vertex $v \in G''$. Then to finish the proof it is enough to show
the existence of an acceptable colouring of the $G''$ from the set
of lists $L''(v)$.

By definition, we have that $\chi(G'')=\chi(G')-r=|{\cal C}|=C$
and $G''$ is a complete $C$ partite graph. The number of vertices
of $G''$ satisfies
$$|V(G'')|=|V(G')|-r<|V(G')|<(2-\epsilon)\chi(G')=(1+o(1))(2-\epsilon)\chi(G'').$$
So by choosing $\delta=\epsilon/2$ we obtain that
$|V(G'')|<(2-\delta)\chi(G'')$. This together with above
discussion implies that  $G''$ satisfies all the condition
($1$--$5$) of the next lemma. This lemma guarantees the existence
of an acceptable colouring of $G''$ from the set of lists $L''(v)$
and completes the proof of Theorem \ref{maint}.

\begin{lemma}
\label{bruce} Let $0<\delta<1$ be a constant and let $C,S,k,m,t$
and $n$ be integers with $m> 6/\delta$, $C=km$, $S=kt$ and $n <
(2-\delta)C$. Suppose, in addition, that $m$ is fixed and $n$ (and
hence $C$) is a sufficiently large function of $m$. Let $G$ be a
complete $C$-partite graph on $n$ vertices and let $L(v)$ be the
set of lists of colours of size $C$ one for each vertex $v$ of $G$
such that the following holds.
\begin{enumerate}
\item
$\cap_{v\in U}L(v)=\emptyset$ for any part $U$ of size bigger than
one in the partition of $G$.
\item
$G$ contains a set of vertices $\cal S$ of size $S$ such that the
vertices in $\cal S$ form parts of size one in the partition of
$G$. The set $\cal S$ is partitioned into $k$ parts ${\cal S}_1,
\ldots, {\cal S}_k$ each of size $t$.
\item
$G$ contains no parts of size two and at most $3k$ singleton parts
which do not belong to $\cal S$.
\item
There exist a set of colours $\cal C$ of size $C$ and its
partition ${\cal C}_1, \ldots, {\cal C}_k$ into $k$ sets of size
$m$. Such that ${\cal C}_i \subset L(s)$ for every vertex $s \in
{\cal S}_i$. In particular, for any subset of ${\cal C}_i$ of size
$t$ there exist an acceptable colouring of the vertices of ${\cal
S}_i$ which uses the colours in this subset.
\item
The total number of colours in the union of all the lists $L(v)$
is less than $n$.
\end{enumerate}
Then there exist an acceptable colouring of $G$ from the lists
$L(v)$.
\end{lemma}

We finish this section with discussion of the proof of Lemma
\ref{bruce}. We postpone all the details to the subsequent
sections of the paper.

\noindent {\bf Proof Overview.}\, The proof proceeds as follows:
\begin{enumerate}
\item[(I)]
We choose a random partition of each ${\cal C}_i$ into two subsets
${\cal A}_i$ of size $t$  and ${\cal B}_i$ of size $m-t$ where
these choices are made uniformly and independently.
\item[(II)]
We use the colours in ${\cal A}_i$ to colour the vertices of
${\cal S}_i$ which is possible by Condition 4 of the lemma.
\item[(III)]
We choose a (random) bijection between ${\cal B}=\cup_{i=1}^k
{\cal B}_i$ and the parts of $G$ not in $\cal S$ in such a way
that, for each part $U$ not in $\cal S$, $U$ is equally likely to
correspond to each  colour $c \in {\cal B}$. We denote by $c_U$ be
the colour corresponding to $U$.
\item[(IV)]
For each part $U$ not in ${\cal S}$ we colour every vertex $v$ of
$U$ for which $c_U \in L(v)$ with the colour $c_U$.
\item[(V)]
We match the set $V^\prime$ of vertices not yet assigned a colour
with the set of colours not yet used (i.e those colours not in
${\cal C}$) so that every vertex is matched with a colour on its
list. We colour each vertex of $V^\prime$ with  the colour with
which it is matched.
\end{enumerate}

If we successfully complete this five step process, we have an
acceptable colouring as every colour not in ${\cal B}$ appears on
at most one vertex, and every colour in ${\cal B}$ appears only on
a subset of some  part, and hence on an independent set of $G$.

To prove that we can find the colouring in this fashion, we need
to describe and analyze our method for choosing the random
bijection between the parts and the colours in ${\cal C}$ made in
Steps I--IV, in order to show that (with positive probability) we
can complete the colouring by finding the desired matching in Step
V.

A key tool will be Hall's Theorem which states that in a bipartite
graph with bipartition $(A,B)$, we can find a matching $M$ such
that every vertex of $A$ is the endpoint of an edge of $M$
provided there is no subset $X$ of $A$ such that setting $N(X)
=\cup \{N(x) | x \in X\}$ we have $|N(X)|<|X|$.

We remark that although Steps I--III are presented as though they
are separate processes performed sequentially, in the more
complicated case of our analysis we will need to interleave these
processes by first choosing some of the ${\cal B}_i$, then
choosing the parts with which these colours will be matched, and
finally completing Step I and then Step III.

To determine if we can find the desired matching in Step V, we
will need to examine the sets $L^\prime(v)=L(v)-{\cal C}$ for the
vertices of $V^\prime$. Let $H$ be  a bipartite graph with
bipartition $\big(V^\prime,\, \cup_{v}L^\prime(v)\big)$ and an
edge between $c$ and $v$ precisely when $c \in L^\prime(v)$. For
each vertex $v$, we let the weight of $v$, denoted $w(v)$, be $1
\over |L^\prime(v)|$. For any set $S$ of vertices we use $W(S)$ to
denote the sum of the weights of the vertices in $S$.

This definition of weight is motivated by the following immediate
consequence  of Hall's Theorem:

\begin{observation}
\label{wobs} If we cannot find the desired matching in Step 5 then
there exist a subset $X$ of $V^\prime$ such that $W(X)>1$. In this
case $W(V^\prime)>1$ as well.
\end{observation}

\noindent {\bf Proof.}\, By Hall's Theorem there exist a subset
$X$ of $V^\prime$ such that $|N(X)|<|X|$. Then, we obtain that
$$
~W(V^\prime) \geq  W(X)=\sum_{x \in X} {1 \over |L^\prime(v)|}=
\sum_{x \in X} {1 \over |N(x) |} \ge \sum_{x \in X} {1 \over |N(X)
|} \ge {|X| \over |N(X)|}>1. ~~~~\Box$$

Thus an analysis of the random parameter $W(V^\prime)$ will be
crucial to the proof of the lemma. In the next section, by
computing the expected value of the parameter, we show that the
lemma holds if $n \le C+S$. In later sections, we complete the
proof using a more complicated analysis along  the same lines.

\section{The expected value of {\boldmath $W(V^\prime)$}}\label{section 3}
\setzero\vskip-5mm \hspace{5mm}

For each part $U$ which is not in ${\cal S}$, our choices in Steps
I and III guarantee that each colour of ${\cal C}$ is equally
likely to be $c_U$. Thus, for each vertex $v$ in such a part, the
probability that $v$ is in  $V^\prime$, i.e., $c_U \not \in L(v)$,
is $1-{|L(v) \cap {\cal C}| \over C}$. Since $|L(v)|=C$, this is
$|L^\prime(v)| \over C$. So, we have:
\begin{equation}
\label{expectation} {\bf E}\big(W(V^\prime)\big)=\sum_{v \in
V-{\cal S}} w(v){\bf Pr}(v \in V^\prime) =\sum_{v \in V-{\cal
S}}\frac{1}{L^\prime(v)} \frac{L^\prime(v)}{C}={n-S \over C}.
\end{equation}
So, if $n \le S+C$, then this expected value is less than or equal
to one. Since the probability that a random variable exceeds its
expected value is less than one,  this implies that we can make
the choices in Steps I--IV so that $W(V^\prime) \le 1$, and hence
by Observation \ref{wobs}, the desired matching can be found in
Step V.

Analyzing the behaviour of the (random) weights of various subsets
of $V^\prime$ will allow us extend our proof technique to handle
larger values of $n$. In doing so, the following definitions and
observations will prove useful.

We let $A$ be the number of non-singleton parts. Then by Condition
3 of the lemma, $A$ is at least $C-S-3k$. Since each non singleton
colour class has at least three vertices and the total number of
classes is $C$, we obtain $(2-\delta)C \ge n \ge C+2A$, i.e., $A
\le {(1-\delta)C \over 2}$. On the other hand, the analysis above
shows that we can assume that $n >S+C$ and hence that $S  \le
(1-\delta) C$. Thus, $A \ge \delta C -3k=\delta C -\frac{3}{m}C
\ge {\delta C \over 2}$. Both these bounds on $A$ will be useful
in our analysis. Note also that
$$m-t=\frac{C-S}{k} \geq \frac{C-(1-\delta) C}{k}=\delta\frac{C}{k}=\delta
m.$$

\section{Completing the proof: the idea}\label{section 4}
\setzero\vskip-5mm \hspace{5mm}

Let $H$ be a bipartite graph with bipartition $\big(V^\prime,\,
\cup_{v}L^\prime(v)\big)$ and an edge between $c$ and $v$
precisely when $c \in L^\prime(v)$. Our first step will be to
check Hall's criterion for a fixed subset of colours $K$ in
$\cup_{v}L^\prime(v)$ and show that the expected number of
vertices in $\{ v| v \in V^\prime, L^\prime(v) \subseteq K\}$ is
less than $|K|$.

To begin, we note that for any such $v$, $w(v) \geq {1 \over |K|
}$. Therefore, defining the set $S_K$ to be $S_K=\{ v| v \in
V^\prime, L^\prime(v) \subseteq K\}$, we have that
\begin{eqnarray*}
{\bf E}\big(|S_K|\big)&=&\sum_{v \in V-{\cal S},\, L^\prime(v)
\subseteq K}{\bf Pr}(v \in V^\prime) \hspace{0.2cm} \leq
\hspace{0.2cm}|K| \sum_{v \in V-{\cal S},\, L^\prime(v) \subseteq
K}
w(v){\bf Pr}(v \in V^\prime)\\
 &\leq& |K| \sum_{v \in V-{\cal S}} w(v){\bf Pr}(v \in
V^\prime)=|K|{\bf E}\big(W(V^\prime)\big)=|K|\frac{n-S}{C}
\end{eqnarray*}
which,  since $n \le (2-\delta)C \le S+A+3k+(1-\delta)C$, is at
most $|K|(1+{A+3k \over C}-\delta )$. On the other hand, this
estimate is not good enough to guarantee Hall's criterion, since
it still can be greater than $|K|$.

To improve on this bound, we use the fact that no colour $c$
appears on the list of all the  vertices of any non-singleton part
of $G$. Note that $k=C/m$, $m>6/\delta$, the number of
non-singleton parts is $A$ and the total number of vertices is at
most $n \le (2-\delta)C \le S+A+3k+(1-\delta)C$. This altogether
implies that for every $c$,
\begin{eqnarray}
\label{weight-color} {\bf E}\Big(W\big(V^\prime \cap \{v | c \in
L^\prime(v) \}\big)\Big) &=& {\bf E}\big(W(V^\prime
)\big)\hspace{0.1cm}- \sum_{v \in V-{\cal S},\, c\not \in
L^\prime(v)}
w(v){\bf Pr}(v \in  V^\prime) \nonumber\\
&=&{\bf E}\big(W(V^\prime )\big)\hspace{0.1cm}- \sum_{v \in
V-{\cal S},\, c\not \in L^\prime(v)}\frac{1}{C}
\le \frac{n-S}{C} -{A \over C}\nonumber\\
&=&\frac{n-S-A}{C} \le
\frac{(1-\delta)C+3k}{C}\nonumber\\
&=&1-\delta+\frac{3}{m} \le 1-{\delta \over 2}.
\end{eqnarray}
Applying this fact for the $c$ in $K$ allows us to improve our
bound on ${\bf E}\big(|S_K|\big)$. Specifically, we note that
summing this bound over all the colours $c$ in $K$
$${\bf E}\Big(\sum_{c \in K} W\big(V^\prime \cap \{v | c \in L^\prime(v)
\}\big)\Big)=\sum_{c \in K}{\bf E}\Big(W\big(V^\prime \cap \{v | c
\in L^\prime(v) \}\big)\Big) \le \left(1 -{\delta \over
2}\right)|K|.$$ Now, each vertex $v$ of $S_K$ contributes
$w(v)=1/|L^\prime(v)|$ to exactly $|L^\prime(v)|$ terms in the
first  sum in this equation, so its total contribution to the sum
is 1. I.e., we have:
$${\bf E}\big(|S_K|\big) \le \left(1 -{\delta \over 2}\right)|K|.$$

So, we don't expect any particular set $K$ of colours to provide
an obstruction to finding the desired matching in the bipartite
graph $H$ in Step V. However, we need to handle all the $K$ at
once. In order to do so, we would like to prove that for each $K$,
the size of $S_K$ is highly concentrated around its expected value
and hence is greater than $|K|$ only with exponentially small
probability. As above, rather than focusing on all the $K$ we
actually consider, for each colour c,
 the weight of the subset $V_c^\prime$ of
$V^\prime$ consisting of those $v$ with $c$ on $L^\prime(v)$.
There are  two major difficulties which complicate our approach.
\begin{itemize}
\item
some of the parts $U$ can be very large making it impossible for
us to obtain the desired concentration results directly (e.g.,
there could be a part of size exceeding $n \over 3$).
\item
If $L^\prime(v)$ is very small then $w(v)=1/|L^\prime(v)|$ is
large and putting $v$ into $V^\prime$ can have a significant
effect on the weight of the various $V_c^\prime$. This makes
proving a concentration result directly impossible.
\end{itemize}

In order to deal with these problems, we proceed as follows:
\begin{enumerate}
\item[(A)]
We colour the ``big'' parts first, ignoring concentration in our
computation and focusing only on the expected weight of the subset
of $V^\prime$ intersecting the big parts. We note that by
considering the expected overall weight and not focusing on a
specific $V_c^\prime$, we only lose a factor of $1 \over C$ per
part. We will define big parts so that there are $o(1)$ of them,
and hence the total loss will not be significant.
\item[(B)]
We treat $v$ with $|L^\prime(v)|$ small  separately using an
expected value argument to bound the weight of the vertices in
this set.
\end{enumerate}

\section{Completing the proof: the details}\label{section 5}
\setzero\vskip-5mm \hspace{5mm}

In this section we will complete the proof of Lemma \ref{bruce}
using the ideas which have already been discussed above. We choose
an integer $b$ so that
$${\delta^2C \over 40 }\le b(m-t) \le {\delta^2C \over 20}$$
which is possible because $m \le {\delta^2 C \over 40}$ (this
holds, since $m$ and $\delta$ are fixed but $C$ tends to infinity)
and $m-t>0$ (in fact it exceeds $\delta m$ as we remarked at the
end of Section \ref{section 3}). We call the largest $b(m-t)$
parts in our partition of $G$ {\it big}, and the others {\it
small}. Let $Big$ be the union of the vertex sets of the big
parts. We will need the following lemma.

\begin{lemma}
\label{Slem} Every small non-singleton partition class contains at
least two $v$ which satisfy:
$$|L^\prime(v)| > {\delta^3 \over 80}C.$$
\end{lemma}

\noindent{\bf Proof.} Let $U$ be a small non-singleton colour
class. We already mentioned that every colour of $\cal C$ is
missed by a vertex of $U$ so $\sum_{v \in U} |L^\prime(v)| =
\sum_{v \in U} |{\cal C}-L(v)| \ge |{\cal C}|=C$. Now, since there
are less than $n< (2-\delta)C$ colours in total, every $L(v)$ must
contain at least $\delta C$ colours in ${\cal C}$ and so the
largest $L^\prime(v)=L(v)-{\cal C}$ in $U$ has at most
$(1-\delta)C$ elements. Thus, the sum of $|L^\prime(v)|$ over the
remaining vertices of $U$ is at least $\delta C$.

Since there are at least $\delta^2 C \over 40$ big colour classes,
the largest small colour class has at most ${40 n \over \delta^2
C} < {80 \over \delta^2}$ vertices. So, the second largest
$L^\prime(v)$ has size at least $\delta C \cdot ({80 \over
\delta^2})^{-1}$.  This is the desired result. \proofend

With this auxiliary result in hand, we can now complete the proof.
We proceed as follows: \vskip0.2cm

\noindent{\bf First Process:} We randomly choose $b$ of the ${\cal
C}_i$ and a partition of each of these into subsets ${\cal A}_i$
of size $t$ and ${\cal B}_i$ of size $m-t$ where these choices are
all made independently and uniformly. We then choose a uniformly
random bijection between the $b(m-t)$ colours in the union of
these ${\cal B}_i$ and the big parts. \vskip0.2cm

\noindent{\bf Second Process:} We chose a partition of each
remaining ${\cal C}_i$ into ${\cal A}_i$ and ${\cal B}_i$ where
again these choices are uniform, independent, and independent of
all the earlier choices. We then choose a uniformly random
bijection between the colours in these ${\cal B}_i$ and the small
parts. \vskip0.2cm

Denote by $c_U$ the colour which is assigned by the above
bijection to the partition class $U$. Use the colours in ${\cal
A}_i$ to colour the vertices of ${\cal S}_i$ and for each part $U$
not in ${\cal S}$ colour every vertex $v$ of $U$ for which $c_U
\in L(v)$ with the colour $c_U$. Let $V^\prime$ be a set of
vertices not yet assigned a colour. We set
$V^{\prime\prime}=V^\prime-Big$ and $V^{\prime\prime\prime}
=V^\prime \cap Big$.

Note that $V^{\prime\prime\prime}$ is determined by our choices in
the first process. So, using a computation similar to that in
(\ref{expectation}) we obtain
$${\bf E}\big(W(V^{\prime\prime\prime})\big)=
\sum_{v \in Big} w(v){\bf Pr}(v \in V^{\prime\prime\prime})
=\sum_{v \in  Big}\frac{1}{L^\prime(v)} \frac{L^\prime(v)}{C}=
{|Big| \over C}.$$ Furthermore, by the definition of expectation,
there exist at least one set of choices for the first process such
that $W(V^{\prime\prime\prime})\le {\bf
E}(W(V^{\prime\prime\prime}))={|Big| \over C}$. We condition on
any such set of choices which ensures that this inequality holds.
We use ${\bf CP}$ and ${\bf CE}$ for the conditional probability
of an event  and conditional expectation of a variable for the
second process, given this set of choices.

Let ${\cal C}^\prime$ be the union of the set of colours in the
${\cal C}_i$ which were chosen in the first process. At the end of
Section 3 we proved that $m-t$ is at least $\delta m$. Therefore
$|{\cal C}^\prime| =mb= {m \over m-t}b(m-t) \leq {\delta}^{-1}
\cdot {\delta^2 C \over 20}={\delta C \over 20}$. Hence, we have
that for every $v$ in a small part which is not in ${\cal S}$,
\begin{eqnarray*}
{\bf CP}(v \in V^{\prime\prime}) &\leq& {|L^\prime(v)| \over
C-|C^\prime|} \le {|L^\prime(v)| \over C}\left(1+ {|C^\prime|
\over C-|C^\prime|}\right)
\leq {|L^\prime(v)| \over C}\left(1+\frac{\delta/20}{1-\delta/20}\right)\\
&\le& \left(1+ {\delta \over 10}\right){|L^\prime(v)| \over C}=
\left(1+ {\delta \over 10}\right) {\bf Pr}(v \in
V^{\prime\prime}).
\end{eqnarray*}

\noindent Clearly, this implies that for every subset $X$ of the
set of vertices $V-{\cal S}-Big$ we have
$${\bf CE}\big(W(V^{\prime\prime} \cap X)\big) \leq \left(1+ {\delta \over
10}\right) {\bf E}\big(W(V^{\prime\prime} \cap X)\big).$$ In
particular, for every colour $c$
\begin{equation}
\label{V''-color} {\bf CE}\Big(W\big(V^{\prime\prime} \cap \big\{
v \,|\,  c \in L^\prime(v)\big\}\big)\Big) \le \left(1+ {\delta
\over 10}\right){\bf E}\Big( W\big(V^{\prime\prime} \cap \big\{ v
\,|\,  c \in L^\prime(v)\big\}\big )\Big)
\end{equation}
and also
\begin{equation}
\label{small-lists1} {\bf CE}\left(W\Big(V^{\prime\prime} \cap
\Big\{ v \, \big|\,|L^\prime(v)| < {n \over \sqrt{\log n}}
\Big\}\Big )\right)\leq \hspace{3.9cm}
\end{equation}
$$\hspace{4.5cm}
\le \left(1+ {\delta \over 10}\right) {\bf E}\left(
W\Big(V^{\prime\prime} \cap \Big\{ v \,\big|\,|L^\prime(v)| < {n
\over \sqrt{\log n}} \Big\}\Big )\right).$$

Before we proceed with the proof, we need the following lemma.
\begin{lemma}
\label{large-deviation} For every color $c$ the probability that
$$W\left(V^{\prime\prime}
\cap \Big\{ v \,\big |\,  c \in L^\prime(v), |L^\prime(v)| \ge {n
\over \log n}\Big\}\right)\hspace{5cm}$$
$$\hspace{3.9cm}
> {\bf CE}\left(W\Big(V^{\prime\prime}
\cap \Big\{ v \, \big| \, c \in L^\prime(v), |L^\prime(v)| \ge {n
\over \log n}\Big\}\Big)\right)+\frac{\delta}{20}$$ is
$o(n^{-1})$.
\end{lemma}

\noindent {\bf Proof.}\, To prove the lemma we need the following
variant of a standard large deviation inequality for martingales.
Since the proof of this inequality is essentially the same as
other proofs which already appeared in the literature (see, e.g.,
Section 3 of the survey \cite{Mc}), we will omit it here.

Given a finite set $\{1, 2, \ldots, r\}$, let $S_r$ denotes the
set of all $r!$ permutations or linear orders on this set. Let
${\bf X}=(X_1, \ldots, X_l)$ be a family of independent random
variables, where the random variable $X_j$ takes values in a
finite set $\Omega_j$. Thus ${\bf X}$ takes values in the set
$\Omega=\prod_j \Omega_j$. Let $\pi \in S_r$ be a random
permutation independent from ${\bf X}$. Suppose that the
non-negative real-valued function $h: \Omega \times S_r
\rightarrow {\mathbb R}$ satisfies the following two conditions
for every $({\bf x}, \pi)$.
\begin{itemize}
\item
For every $j$, changing the value of a coordinate $x_j$ can change
the value of $h({\bf x}, \pi)$ by at most $d$.
\item
Swapping any two elements in permutation $\pi$ can change the
value of $h({\bf x}, \pi)$ by at most $d$.
\end{itemize}
Denote by ${\bf E}h$ the expected value of $h$. Then for every $t
\geq 0$ we have that
$${\bf Pr}\Big(|h-{\bf E}h|>t\Big) \leq
e^{-\Omega\big(\frac{t^2}{(r+l)d^2}\big)}.$$

Now fix a color $c$ and define the function $h$ to be
$$h({\bf x}, \pi)= W\left(V^{\prime\prime}
\cap \Big\{ v \,\big |\,  c \in L^\prime(v), |L^\prime(v)| \ge {n
\over \log n}\Big\}\right),$$ where $({\bf x}, \pi)$ corresponds
to the set of random choices for the second process. More
precisely, $x_i$ is a random partition of the set ${\cal C}_i$
into subsets ${\cal A}_i$ and ${\cal B}_i$ and $\pi$ is a random
bijection between the colors in these ${\cal B}_i$ and the small
parts of $G$. Since we can fix one canonical ordering of these
small parts we can assume that $\pi$ is just a random permutation
of the set of colors which is, by definition, independent from the
variables $x_i$.

Next, note that changing the outcome of the variable $x_i$, i.e.,
changing one particular ${\cal B}_i$ can only affect vertices in
at most $m-t$ small parts of $G$. As we already mentioned in the
proof of Lemma \ref{Slem}, each small part contains at most
$80/\delta^2$ vertices. Since we considering only vertices $v$
satisfying $|L^\prime(v)|\geq \frac{n}{\log n}$, the weight of
such a vertex is at most $w(v)=1/|L^\prime(v)| \leq \frac{\log
n}{n}$. Therefore, changing outcome of one $x_i$ can change the
value of $h$ by at most $(m-t)\frac{80}{\delta^2}\frac{\log n}{n}=
O\big(\frac{\log n}{n}\big)=d$. Similarly swapping any two colors
in $\pi$ can affect only vertices in two small parts of $G$. So
again this can only change $h$ by at most $d=O\big(\frac{\log
n}{n}\big)$.

Since the total number of random variables $x_i$ and also the
length of permutation $\pi$ are bounded by $n$ we have that in our
case $(r+l)d^2 \leq O\big(n(\frac{\log n}{n})^2\big)=
O\big(\frac{\log^2 n}{n}\big)$. Therefore it follows form the
above large deviation inequality that
$${\bf Pr}\Big(h-{\bf E}h >t=\frac{\delta}{20}\Big) \leq
e^{-\Omega\big(\frac{t^2}{(r+l)d^2}\big)}=
e^{-\Omega\big(\frac{(\delta/20)^2}{O(\log^2 n/n)}\big)}=
e^{-\Omega\big(\frac{n}{\log^2 n}\big)}=o(n^{-1}).$$ This
completes the proof of the lemma. \hfill $\Box$

Now, using the fact that the total number of colors is at most
$n$, we deduce from this lemma that with probability $1-o(1)$ the
following holds for every color $c$
\begin{equation}
\label{concentration} W\left(V^{\prime\prime} \cap \Big\{ v \,\big
|\,  c \in L^\prime(v), |L^\prime(v)| \ge {n \over \log n}
\Big\}\right)  \hspace{4.5cm}
\end{equation}
$$ \hspace{3cm } \leq
 {\bf CE}\left(W\Big(V^{\prime\prime}
\cap \Big\{ v \, \big| \, c \in L^\prime(v), |L^\prime(v)| \ge {n
\over \log n} \Big\}\Big)\right)+\frac{\delta}{20}.$$ In addition,
we also want to satisfy the following inequality:
\begin{equation}
\label{small-lists2} W\left(V^{\prime\prime} \cap \Big\{ v \,
\big| \,|L^\prime(v)| < {n \over \sqrt{\log n}} \Big\}\right )
 \hspace{5cm}
\end{equation}
$$\hspace{3.5cm}
\leq \left(1+ {\delta \over 10}\right){\bf
CE}\left(W\Big(V^{\prime\prime} \cap \Big\{ v  \,\big|
\,|L^\prime(v)| < {n \over \sqrt{\log n}} \Big\}\Big)\right).$$
Since the probability that this last inequality fails is at most
$\frac{1}{1+\delta/10}<1-o(1)$, there does indeed exist a set of
random choices for the second process which satisfies
simultaneously (\ref{concentration}) and (\ref{small-lists2}).

Fix any such set of choices. Then, combining the inequalities
(\ref{small-lists1}) and (\ref{small-lists2}) together with the
facts that $W(V^{\prime\prime\prime})\leq \frac{|Big|}{C}$ and
$V^\prime= V^{\prime\prime}\cup V^{\prime\prime\prime}$ we obtain
that
$$W\left(V^\prime \cap \Big\{v \, \big| \, |L^\prime(v)| < {n \over
\sqrt{\log n}}\Big\} \right)\leq W\left(V^{\prime\prime} \cap
\Big\{v \, \big| \, |L^\prime(v)| < {n \over \sqrt{\log n}}\Big\}
\right)+ W\big(V^{\prime\prime\prime}\big)$$
\begin{eqnarray}
\label{small-sets} \hspace{1.7 cm} \leq \left(1+ {\delta \over
10}\right)^2 {\bf E}\left( W\Big(V^{\prime\prime} \cap \Big\{ v
\,\big|\,|L^\prime(v)| < {n \over \sqrt{\log n}} \Big\}\Big
)\right)+\frac{|Big|}{C}.
\end{eqnarray}
Note that, by Lemma \ref{Slem}, every small non-singleton
partition class contains at least two vertices $v$ such that
$|L^\prime(v)|>\frac{\delta^3}{80}C=\Omega(n)> \frac{n}{\sqrt{\log
n}}$. Since the number of small non-singleton partition classes is
at least $A-\frac{\delta^2}{20}C$ we obtain that
$$\hspace{-1cm}\left|(V-{\cal S}-Big) \cap \Big\{ v \,\big|\,|L^\prime(v)|
< {n \over \sqrt{\log n}} \Big\}\right| \hspace{0.2cm}\leq
\hspace{0.2cm}n-S-|Big|- 2\Big(A-\frac{\delta^2}{20}C\Big)$$
\begin{eqnarray*}
\hspace{5cm}&\leq&
(2-\delta)C-S-|Big|-2A+\frac{\delta^2}{10}C\\
&=&
(1-\delta)C-|Big|+(C-S-A)-A+\frac{\delta^2}{10}C\\
&\leq&
(1-\delta)C-|Big|+3k-A+\frac{\delta^2}{10}C \\
&\leq& \left(1-\frac{4}{5}\delta\right)C-|Big|.
\end{eqnarray*}
Here, in the last inequality we used that
$A>\frac{\delta}{2}C>\frac{3}{m}C=3k$ and $\delta^2 \leq \delta$.
Note that a similar computation as in (\ref{expectation}) shows
that for any subset $Y \subseteq V-{\cal S}$ the expectation ${\bf
E}\big(W(V^\prime \cap Y)\big)=\frac{|Y|}{C}$. In particular, for
$Y=\big(V-{\cal S}-Big\big) \cap \big\{ v \,\big|\,|L^\prime(v)| <
{n \over \sqrt{\log n}} \big\}$ we obtain
\begin{eqnarray*}
{\bf E}\left( W\Big(V^{\prime\prime} \cap \Big\{ v
\,\big|\,|L^\prime(v)| < {n \over \sqrt{\log n}} \Big\}\Big
)\right) \hspace{-0.1cm}&=&\hspace{-0.1cm} \frac{\big|(V-{\cal
S}-Big) \cap \{ v \,|\,|L^\prime(v)|
< {n \over \sqrt{\log n}} \}\big| }{C} \\
\hspace{-0.1cm}&\leq&\hspace{-0.1cm}
1-\frac{4}{5}\delta-\frac{|Big|}{C}.
\end{eqnarray*}
Combining this inequality with (\ref{small-sets}) we have
\begin{eqnarray}
\label{small-sets1} W\left(V^\prime \cap \Big\{v \, \big| \,
|L^\prime(v)| < {n \over \sqrt{\log n}}\Big\} \right) &\leq&
\left(1+ {\delta \over 10}\right)^2
\left(1-\frac{4}{5}\delta-\frac{|Big|}{C}\right)+\frac{|Big|}{C}
\nonumber\\
&\leq&
\left(1+\frac{\delta}{4}\right)\left(1-\frac{4}{5}\delta\right)
\le 1-{\delta \over 2}.
\end{eqnarray}
This completes our analysis of the weight of vertices with short
lists. We now consider the remaining vertices.

As we already mentioned, for every color $c$ and every
non-singleton part of $G$ there is at least one vertex $v$ in this
part such that $c \not \in L(v)$. Since there are at least
$A-\frac{\delta^2}{20}C$ small non-singleton parts, a similar
computations as in (\ref{weight-color}) shows for every color $c$
that
\begin{eqnarray*}
{\bf E}\bigg(W\big(V^{\prime\prime} \cap \big\{v \, \big|\, c \in
L^\prime(v)\big\}\Big)
\bigg)&=&\frac{\big|(V-{\cal S}-Big) \cap \{v \, |\, c \in L^\prime(v)\}\big|}{C}\\
&\leq& \frac{n-S-|Big|-(A-\frac{\delta^2}{20}C)}{C} \\
&\leq& \frac{(2-\delta)C-S-|Big|-A}{C}+\frac{\delta^2}{20}\\
&=&
(1-\delta)+\frac{C-S-A}{C}-\frac{|Big|}{C}+\frac{\delta^2}{20}\\
&\leq&(1-\delta) +\frac{3k}{C} -\frac{|Big|}{C} +\frac{\delta^2}{20}\\
&=&
(1-\delta)+\frac{3}{m}-\frac{|Big|}{C}+ \frac{\delta^2}{20}\\
\end{eqnarray*}
$$ \hspace*{5cm} \leq  1-\delta+\frac{\delta}{2}+\frac{\delta^2}{20}-\frac{|Big|}{C}\leq 1-\frac{2}{5}\delta
-\frac{|Big|}{C}. $$
Combining this inequality with (\ref{V''-color}) and (\ref{concentration}) and using the fact
that $V^\prime=V^{\prime\prime} \cup V^{\prime\prime\prime}$ we will have that for every color $c$
$$\hspace{-5.9cm}W\left( V^\prime \cap
\Big\{ v \big|  c \in L^\prime(v), |L^\prime(v)| \ge {n \over \log
n} \Big\} \right) $$
\begin{eqnarray}
\label{all-colors} \hspace{1.8cm} &\leq & W \left(
V^{\prime\prime} \cap \Big\{ v \big|  c \in L^\prime(v),
|L^\prime(v)| \ge {n \over \log n} \Big\} \right)
+W\big(V^{\prime\prime\prime}\big)\nonumber\\
&\leq& {\bf CE}\left( W \Big(V^{\prime\prime} \cap \Big\{ v \big|
c \in L^\prime(v), |L^\prime(v)| \ge {n \over \log n}
\Big\}\Big)\right)
+\frac{\delta}{20} +\frac{|Big|}{C}\nonumber\\
&\leq& {\bf CE}\bigg( W \left( V^{\prime\prime} \cap \Big\{ v
\big|  c \in L^\prime(v) \Big\}\right)\bigg)
+\frac{\delta}{20} +\frac{|Big|}{C}\nonumber\\
&\leq& \left(1+ {\delta \over 10}\right){\bf E}\bigg( W \left(
V^{\prime\prime} \cap \Big\{ v \big|  c \in L^\prime(v)
\Big\}\right)\bigg)
+\frac{\delta}{20} +\frac{|Big|}{C}\nonumber\\
&\leq& \left(1+ {\delta \over 10}\right)\left( 1-\frac{2}{5}\delta
-\frac{|Big|}{C}\right)
+\frac{\delta}{20} +\frac{|Big|}{C}\nonumber\\
&\le& \left(1+ {\delta \over 10}\right)\left(
1-\frac{2}{5}\delta\right)+\frac{\delta}{20} \leq 1-{\delta \over
4}.
\end{eqnarray}

Recall that $H$ is  a bipartite graph with bipartition
$\big(V^\prime,\, \cup_{v}L^\prime(v)\big)$ and an edge between
$c$ and $v$ precisely when $c \in L^\prime(v)$. Let $K$ be any
subset of colours in $\cup_{v}L^\prime(v)$ and denote by $S_K=\{
v| v \in V^\prime, L^\prime(v) \subseteq K\}$. We complete the
proof of the lemma by showing that the graph $H$ satisfies Hall's
condition, i.e., $|S_K| \leq |K|$ for every set $S_K$. Then in
Step V we can match all uncoloured vertices in $V^\prime$ with the
set of colours yet not used and and produce an acceptable coloring
of $G$.

First, note that any set $K$ of fewer than $n \over \sqrt{\log n}$
colours cannot be an obstruction to the existence of the desired
matching. Indeed, if $|S_K|> |K|$, then by Observation \ref{wobs}
we have that $W(S_K)>1$. On the other hand, for every vertex $v
\in S_K$ the size of $L^\prime(v)$ is at most $|K|<{n \over
\sqrt{\log n}}$. Therefore we obtain a contradiction, since by
(\ref{small-sets1})
$$W\big(S_K\big) \leq W\left(V^\prime \cap \Big\{v \, \big| \, |L^\prime(v)| <
{n \over \sqrt{\log n}}\Big\} \right ) <1-{\delta \over 2}.$$
\noindent Turning to larger $K$, we note next that the inequality
(\ref{small-sets1}) yields:
$$\left|V^\prime\hspace{-0.02cm} \cap\hspace{-0.02cm} \Big\{v\, \big|\,
|L^\prime(v)| \leq {n \over {\log n}}\Big\}\right| \le
W\hspace{-0.05cm}\left(\hspace{-0.03cm}V^\prime
\hspace{-0.02cm}\cap \hspace{-0.02cm}\Big\{v \, \big|\,
|L^\prime(v)| \leq {n \over {\log
n}}\Big\}\hspace{-0.03cm}\right)\hspace{-0.06cm} \bigg(\min_{v,
\,|L^\prime(v)| \leq {n \over {\log n}}}
\hspace{-0.06cm}w(v)\hspace{-0.03cm} \bigg)^{-1}$$
\begin{eqnarray}
\label{small-sets2} \hspace{1cm} \leq \left( 1-{\delta \over
2}\right) \bigg(\min_{v,\, |L^\prime(v)| \leq {n \over {\log
n}}}\, \frac{1}{|L^\prime(v)|} \bigg)^{-1} \leq \left( 1-{\delta
\over 2}\right){n \over \log n}<{n \over \log n}.
\end{eqnarray}
Next, observe that the set of inequalities (\ref{all-colors})
imply that for any set of colours $K$
\begin{eqnarray*}
\left|S_K \cap \Big\{v \,\big |\,  |L^\prime(v)| >{n \over \log
n}\Big\}\right|
&=&\sum_{v \in S_K, |L^\prime(v)| >{n \over \log n}} \, w(v)\cdot |L^\prime(v)|\\
&\leq& \sum_{c \in K}\,\,\sum_{\{v\,|\, c\in L^\prime(v),\,
|L^\prime(v)| >{n \over \log
n}\}}\, w(v)\\
&=&\sum_{c \in K} W \left( V^\prime \cap \Big\{ v \,\big|\,  c \in
L^\prime(v), |L^\prime(v)| \ge {n \over \log n} \Big\}
\right)\\
&\le& \left(1-{\delta \over 4}\right)|K|.
\end{eqnarray*}
\noindent This, together with the inequality (\ref{small-sets2})
yields that any set of colours $K$ of size  at least ${n \over
\sqrt{\log n}}$ satisfies
\begin{eqnarray*}
\big|S_K \Big|&=&\left|S_K \cap \Big\{v\,\big |\, |L^\prime(v)|
>{n \over \log n}\Big\}\right|+
\left|S_K \cap \Big\{v \,\big |\, |L^\prime(v)| \leq {n \over \log n}\Big\}\right|\\
&\le& \left(1-{\delta \over 4}\right)|K|+ \left|V^\prime \cap
\Big\{v \,\big |\,
|L^\prime(v)| \leq {n \over \log n}\Big\}\right| \\
&\leq& \left(1-{\delta \over 4}\right)|K|+ {n \over \log n}\, <
\,|K|.
\end{eqnarray*}
Thus we obtain that these larger $K$ also do not violate Hall's
condition and hence the desired matching of Step V does indeed
exist. This completes the proof. \hfill $\Box$

\section{Concluding remarks}\label{section 6}
\setzero\vskip-5mm \hspace{5mm}

In this paper we proved that for every $\epsilon>0$ and for every
sufficiently large graph $G$ of order $n$, the list chromatic
number of $G$ equals its chromatic number, provided $n \leq
(2-\epsilon)\chi(G)$. A more careful analysys of our methods
yields that the value of $\epsilon$ in this result can be made as
small as $O(1/\log^{\eta} n)$ for any constant $0 < \eta<1$.
Nevertheless the conjecture of Ohba remains open for  graphs with
$2\chi(G)$ vertices and it seems one needs new ideas to tackle
this problem. Even to show that there is a constant $N$ such that
$\chi^l(G)=\chi(G)$ for every graph $G$ with at most $2\chi(G)-N$
vertices, would be very interesting.

In conclusion we would like to propose a related problem, which
was motivated by Ohba's conjecture. Let $t$ be an integer and let
$G$ be a graph with at most $t \chi(G)$ vertices. Find the
smallest constant $c_t$ such that for any such a graph $G$ its
list chromatic is bounded by $c_t\chi(G)$. Note that Ohba's
conjecture if true, implies that $c_2=1$. An additional intriguing
question is to determine graphs with $|V(G)| \leq t \chi(G)$ and
for which the ratio $\chi^l(G)/\chi(G)$ is maximal.  Here the case
$t=2$ gives some indication that a complete multi-partite graph
with all parts of size $t$ may have this property.

\label{lastpage}

\end{document}